\documentstyle[amssymb]{amsart}
\input{psfig}
\numberwithin{equation}{section}

\newcommand{\R}{\Bbb{R}}

\newcommand{\tends}{\rightarrow}
\theoremstyle{plain}

\theoremstyle{remark}

\newtheorem{remark}{Remark}

\begin{document}
\title{Picturing Pinchuk's Plane Polynomial Pair}
\author{L. Andrew Campbell}
\address{The Aerospace Corp. \\
M1-102 PO Box 92957 \\
Los Angeles CA 90009}
\email{campbell@@aero.org}
\date{May 9, 1997}
\maketitle
\begin{abstract}
Sergey Pinchuk discovered a class of pairs of real polynomials
in two variables that have a nowhere vanishing Jacobian
determinant and define maps of the real plane to itself
that are not one-to-one. This paper describes the asymptotic
behavior of one specific map in that class.
The level of detail presented permits a good geometric
visualization of the map.
Errors in an earlier description of the image of the map
are corrected (the complement of the image consists of two,
not four, points). Techniques due to Ronen Peretz are
used to verify the description of the asymptotic variety
of the map.
\end{abstract}

\section{Introduction}
The strong real Jacobian conjecture stated 
that every polynomial map from $\R^n$ to $\R^n$ with nowhere vanishing
Jacobian determinant is univalent (one-to-one).
This conjecture was refuted (for $n=2$ and hence all larger $n$) in 1994 by Sergey
Pinchuk, who provided a class of counterexamples
\cite{StrongRealJC}.
One of the counterexamples
is a map
$F=(P,Q): \R^2 \to \R^2$
with $P(x,y)$ and $Q(x,y)$ polynomials of total degree $10$ and $25$, respectively
\cite{Pinchuk.Cx.Email}.
That particular map is the primary focus of this paper.
It can be described as follows.

Let $t = xy-1, h = t(xt+1),
f = ((h+1)/x)(xt+1)^2, P = f+h,
Q=-t^2-6th(h+1)-170fh-91h^2-195fh^2-69h^3-75fh^3-(75/4)h^4$.
Then $F(x,y) = (P(x,y),Q(x,y))$ is a real polynomial
map from $\R^2$ to itself; its Jacobian determinant 
is everywhere positive; and it is not univalent.

This map has been considered elsewhere, in particular in
\cite{PPR2}, where the assertion was made that $F(\R^2)$
consists of all of $\R^2$ except for exactly four points.
There were errors and oversights in the calculations,
and two of the four points cited are actually in the image.
In this paper, the complement of the image is identified as
consisting of the points $(0,0)$ and $(-1,-163/4)$ only.

The asymptotic behavior of the map is studied. In particular,
the asymptotic variety of $F$ (as defined by Ronen Peretz in
\cite{AsympVals}
) is computed using Peretz's technique. 
Denote it by $\text{AV}[F]$.
For the particular $F$ studied here, $\text{AV}[F]$
admits a parameterization by two polynomials of degree two and five
in a single variable.

The {\it asymptotic
flower} of $F$ (new terminology, see
\cite{Flower}
) is the inverse image under $F$ of $\text{AV}[F]$.
Denote it by $\text{AF}[F]$.
By construction, the restriction of $F$ to a mapping
from $\R^2 \setminus \text{AF}[F] $ to $\R^2 \setminus \text{AV}[F]$ is a
proper map. For this particular $F$, it reduces to
homeomorphisms of four
simply connected domains in $\R^2$ , each mapping onto one of the two
simply connected components of $\R^2 \setminus \text{AV}[F]$.
This description provides a good geometric visualization
of the map (and is supplemented by graphics).

\section{Asymptotics of Pinchuk's Map}
Pinchuk's map $F(x,y)=(P(x,y),Q(x,y))$ is most easily
studied by considering the fibers $P=c$ of the map $P$,
because $P$ only has degree 10, whereas $Q$ has degree
25.
The following information and table are excerpted from \cite{PPR2}.
The fiber $P=0$ has five components and $P=-1$ has four components.
In both cases ($c=0$ and $c=-1$) the fibers can be computed and
their components parameterized explicitly without great
difficulty, because the polynomial $P-c$ factorizes simply.
The other fibers
are parameterized by the rational curve
$$x(h) = \frac{ (c-h)(h+1) }{  (c-2h-h^2)^2 }$$
$$ y(h) = \frac{ (c-2h-h^2)^2(c-h-h^2) }{ (c-h)^2 }$$
For a fixed value $c$, the components of the fiber $P=c$ are the images
the map $h \mapsto (x(h),y(h))$ for values of $h$
between successive poles
(which occur when $h=c$ or $c-2h-h^2=0$; no cancellation
occurs as long as $c$ is neither $0$ nor $-1$).
The table below summarizes the data on number of components and
the range of $Q$ for all fibers $P=c$.
$Q$ is always monotone (hence one-to-one) on any component
of a fiber $P=c$, because the Jacobian determinant of $P$ and $Q$
is everywhere nonzero.

\begin{table}[htb]
\label{PinTbl}
\begin{center}
\begin{tabular}{||c||c||}
\hline
$P=c$  & Ranges of Q on the components \\
\hline \hline
$c>0$ & $(+\infty,q-),\; (q-,q+),\; (q+,-\infty),\; (-\infty,+\infty)$ \\
\hline
$c=0$ & $(0,208),\; (-\infty,0),\; (0,+\infty),\; (-\infty,0),\; (208,+\infty)$
\\
\hline
$-1<c<0$ & $(+\infty,q-),\; (q-,-\infty),\; (-\infty,q+),\; (q+,+\infty)$ \\
\hline
$c=-1$ & $(-\infty,-163/4),\; (-\infty,-163/4),\; (-163/4,+\infty),\; (-163/4,+\infty)$ \\
\hline
$c<-1$ & $(+\infty,-\infty),\; (-\infty,+\infty)$ \\
\hline \hline
\multicolumn{2}{|| c ||}{Legend: $(a,b)$ denotes the
open interval from $\min(a,b)$ to $\max(a,b)$;} \\
\multicolumn{2}{|| c ||}
{$q+$ ($q-$) = the value of $Q$ at $h =-1+\sqrt{1+c}$ (resp., $-1-\sqrt{1+c}$);}
 \\
\hline \hline
\end{tabular}
\end{center}
\caption{Ranges of $Q$ on fibers $P=c$ for Pinchuk's map}
\end{table}

\begin{remark}
In \cite{PPR2} there was a typographical error in the formula
for $x(h)$ (the term in the denominator was not squared). The
computations leading to the results in Table 1 used the correct
parameterization (the one shown above). Also, in \cite{PPR2}
one of the points listed as not in the image of $F$ was the
point $(0,208)$. However, a glance at the table shows that
this cannot be correct; the value $208$ lies in $(0,+\infty)$,
which is the range of $Q$ on one the components of the fiber
$P=0$.
\end{remark}

\begin{remark}
In \cite{PPR2}, the parameterization by $x(h),y(h)$ was introduced
without any indication of how it arises. It comes from a straightforward
process of solving the equations that define $P$, first for $x$ and
then for $y$. For example, if $P=c$ then the first step is
$c = f + h$, then $c-h = ((h+1)/x)(xt+1)^2 = ((h+1)/x)(h/t)^2$,
hence $xt^2 = (h+1)h^2/(c-h)$. From the defining equations again,
$t = h - xt^2$, which allows solving for $t$ in terms of $h$,
then for $x = xt^2/t^2$, and finally for $y=(t+1)/x$.
\end{remark}

The finite endpoints of ranges of $Q$ occur precisely due to
components of a fiber along which the $x$ or $y$ component
blows up, but $Q(x,y)$ does not.
Denote $x(h),y(h)$ by $x(c,h),y(c,h)$ to capture the dependence
on $c$. Then one has the following rational identities
$$ P(x(c,h),y(c,h)) = c $$
\begin{multline}
 Q(x(c,h),y(c,h)) =
\frac{1}{4(c-h)^2}
\lbrace 197h^6 +(416-726c)h^5 \\
+ (252 -1684c +825c^2)h^4
+ (-1224c +2040c^2 -300c^3)h^3 \\
+ (1648c^2 -780c^3)h^2
+ (-680c^3)h \rbrace \notag
\end{multline}
The identities can be verified simply by substitution.
For $c \ne 0, -1$ it can be checked that $Q$ blows up
when $h$ tends to $c$ (one of the poles of the parameterization),
but not when $h$ tends to either of the values $-1+\sqrt{1+c}$,
$-1-\sqrt{1+c}$ (which are also poles of the parametrization,
as they are the zeroes of $c-2h-h^2$). The respective values
of $q$ are denoted by $q+$,$q-$. Of course, they depend on $c$.
By definition, the asymptotic variety of a map \cite{AsympVals}
consists of points in the image plane that are limits of the
images of points along a curve that tends to infinity in the
original plane.
By that definition,
for each $c \ne 0,-1$ the points $(c,q+)$ and $(c,q-)$ are in
the asymptotic variety, $\text{AV}[F]$, of the map $F$. These
points can be obtained by simply substituting $c=2h+h^2$ into
the above rational identities for $P$ and $Q$. To make life simple,
$u$ and $v$ will be used as coordinates in the image plane, with
$x$ and $y$ reserved for points in the original plane. Carrying
out the indicated substitution yields the following parameterized
curve in the image plane
$$ u = P = c = 2h + h^2 $$
$$ v = Q = -(1/4)(1736h^3 + 1044h^2 + 1155h^4 + 300h^5)$$
The values of $h$ that lead to $c = 0,-1$ are $h=-1$ ($c=0$),
$h=0$ ($c=0$), and $h=-2$ ($c=0$). The corresponding points
$(u,v)$ arising from the above parameterization are, respectively,
$(-1,-163/4)$,$(0,0)$, and $(0,208)$. From Table 1, these points
all belong to the $\text{AV}[F]$.
So the entire curve lies in $\text{AV}[F]$. 
Using Peretz's technique of standard asymptotic identities, it
will be shown below that this is, in fact, the entire asymptotic
variety $\text{AV}[F]$. Figure 1 is a depiction of the variety.

\begin{figure}[hbt]
\vspace{-0.5in}
\centerline{
\psfig{figure=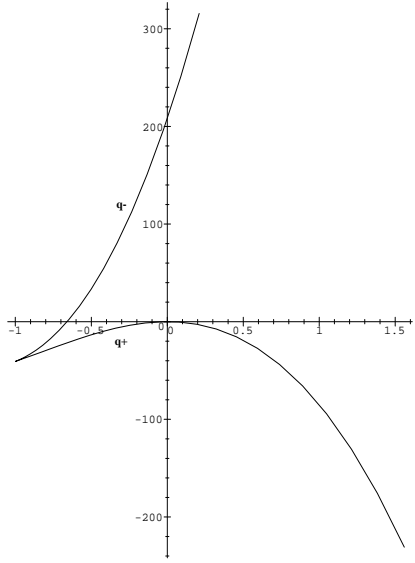,height=4in,width=3in}
}
\caption{The asymptotic variety of Pinchuk's map.}
\label{AVfig}
\end{figure}

The figure illustrates the fact that $q-$ is the larger of the
two values of $Q$ for a given $c$, whereas $q+$ is the smaller
(except at $c=-1$, where they coincide). This can be seen easily
by using a parametrization in terms of $w = h+1$, for which the
upper portion corresponds to $w < 0$ and the lower portion to $w > 0$.

\begin{remark}
In \cite{PPR2} it was claimed that there was a point of the form
$u=c,v=d$ for some $-1 < c < 0$ which was not in the image of $F$,
because the values of $q+$ and $q-$ supposedly coincided at that
point. (That is, $d=q+=q-$, hence $d$ would not lie in any of
the ranges of $Q$ on the fiber $P=c$.)
From the above figure, it is clear that there is no point
where $q+$ or $q-$ coincide for $c > -1$. 
\end{remark}

Table 1 can be rewritten using the fact that $q+ < q-$ for
$c \ne 0,-1$ to put all the intervals $(a,b)$ that are ranges
of $Q$ on fibers of $P$ in canonical form - that is, with $a < b$.
The result is

\begin{table}[htb]
\label{PinTbl2}
\begin{center}
\begin{tabular}{||c||c||}
\hline
$P=c$  & Ranges of Q on the components \\
\hline \hline
$c>0$ & $(-\infty,q+),\; (q+,q-),\; (q-,+\infty),\; (-\infty,+\infty)$ \\
\hline
$c=0$ & $(0,208),\; (-\infty,0),\; (0,+\infty),\; (-\infty,0),\; (208,+\infty)$\\
\hline
$-1<c<0$ & $(-\infty,\; q+),\; (q+,\; +\infty),\; (-\infty,\; q-),\; (q-,\; +\infty)$ \\
\hline$c=-1$ & $(-\infty,-163/4),\; (-\infty,-163/4),\; (-163/4,+\infty),\; (-163/4,+\infty)$ \\
\hline
$c<-1$ & $(-\infty,+\infty),\; (-\infty,+\infty)$ \\
\hline \hline
\end{tabular}
\end{center}
\caption{Ranges of $Q$ on fibers $P=c$ for Pinchuk's map - rewritten}
\end{table}

This clearly shows that the only points omitted from the image of
$F$ are the points $(-1,-163/4)$ and $(0,0)$.

\section{The Peretz Method}
This section uses the techniques described in Ronen Peretz's
paper \cite{AsympVals} to derive conditions that must be
satisfied by any asymptotic values of the polynomial $P$. 
In the next section,
those conditions will be used to show that $\text{AV}[F]$ is exactly the
curve shown in Figure 1.

Observe
first that the highest (total) degree term in $P$ is $x^6y^4$, so $P$
satisfies the Peretz normalization criterion $\deg(P) = \deg_x(P) + \deg_y(P)$.
This implies that $P$ has only $x$ or $y$-finite asymptotic curves.
In fact, $P$ can only have asymptotic curves with $x \tends \pm\infty$
and $y \tends 0$, or vice versa.
To search for $y$-finite asymptotic curves and the corresponding asymptotic
values, follow the steps outlined by Peretz. That is, first write
\begin{equation*}
P(x,y) = P_6x^6 + P_5x^5 + P_4x^4 + P_3x^3 + P_2x^2 + P_1x + P_0
\end{equation*}
where $P_0,\ldots,P_6$ are polynomials in $y$.
This yields
\begin{align}
P_6 &= y^4 \notag \\ P_5 &= -4y^3 \notag \\ P_4 &= 3y^3 + 6y^2 \notag \\
P_3 &= -7y^2 - 4y \notag \\ P_2 &= 3y^2 + 5y + 1 \notag \\
P_1 &= -3y -1 \notag \\ P_0 &= y \notag
\end{align}
Then, assuming that $P$ tends to the (finite) value $C$ along an
asymptotic curve (one that tends to infinity in the domain space), 
write down the Peretz assertions
\begin{align}
&P_6x^6 + P_5x^5 + P_4x^4 + P_3x^3 + P_2x^2 + P_1x + P_0(0) \tends C \notag \\
&P_6x^5 + P_5x^4 + P_4x^3 + P_3x^2 + P_2x + P_1(0) \tends 0 \notag \\
&P_6x^4 + P_5x^3 + P_4x^2 + P_3x + P_2(0) \tends 0 \notag \\
&P_6x^3 + P_5x^2 + P_4x + P_3(0) \tends 0 \notag \\
&P_6x^2 + P_5x + P_4(0) \tends 0 \notag \\
&P_6x + P_5(0) \tends 0 \notag \\
&P_6(0) \tends 0  \notag
\end{align}
These follow from the fact that if a product of two factors tends
to a finite limit and one factor is $x$ (which tends to $\pm\infty$),
then the other factor tends to $0$.

Look, from the bottom up, for the first assertion in which the
term before $\tends$ is not zero. This is the assertion containing $P_2(0)$.
Write the assertion out in full. After judicious factorization it is
\begin{equation*}
(yx)^4 - 4(yx)^3 + (3y+6)(yx)^2 + (-7y - 4)(yx) + 1 \tends 0
\end{equation*}
and since $y \tends 0$, this implies that $xy \tends r$, where
$r$ is a root of
\begin{equation*}
r^4 -4r^3 + 6r^2 -4r + 1 = (r-1)^4
\end{equation*}

\begin{remark}
Justification. First, the power product $xy$ must remain bounded, otherwise the
expression would not approach a finite limit. 
Next, even if the polynomial has multiple distinct roots,
the image of the curve must ultimately remain near a single
root, otherwise the value of the expression would not be ultimately
small (it would not tend to zero because of the trips from one root
to another, necessarily involving points away from the roots).
Similar reasoning applies later, when other power products than $xy$ are considered.
\end{remark}

Thus $xy \tends 1$, which means that $(y-1/x)x \tends 0$. In other
words
\begin{equation*}
y = 1/x + o(1/x)
\end{equation*}
Next denote the error term by $z$, so that $y=1/x+z$. Substitute
$y=1/x+z$ into the Peretz assertions to obtain
the following ones
\begin{align}
3z^2x^2 + x^6z^4 + 2z^2x^3 + 3z^3x^4 + 3zx \tends C \notag \\
z^4x^5 + 3z^3x^3 + 2z^2x^2 + 3z^2x + 6z + 3/x \tends 0 \notag \\
z^4x^4 + 3z^3x^2 + 2z^2x -5z - 4/x \tends 0 \notag \\
x^3z^4 + (-3z^2+z^2(3+3z))x - 2z + 2z(3+3z) + 3z^2 + (3+9z)/x + 3/x^2 \tends 0 \notag \\
z^4x^2 -6z^2 -8z/x -3/x^2 \tends 0 \notag \\
z^4x + 4z^3 + 6z^2/x + 4z/x + 1/x^3 \tends 0 \notag \\
0 \tends 0  \notag
\end{align}
Using the facts that $1/x \tends 0$, $z \tends 0$, and $xz \tends 0$,
these assertions can be immediately simplified to
\begin{align}
x^6z^4 + 2x^3z^2 + 3x^4z^3 &\tends C \notag \\
x^5z^4 &\tends 0 \notag 
\end{align}
plus five additional trivial assertions of the form
$0 \tends 0$.
The fact that $x^5z^4 \tends 0$ imposes the requirement that
$z^4 = o(x^{-5})$. 
It follows that $z = o(|x|^{-5/4})$.
No specific data on the form of the error term is implied.
Finally, the assertion $(x^3z^2)^2 + (2+3xz)(x^3z^2) \tends C$,
together with $xz \tends 0$, means that $x^3z^2$ 
tends to a root $r$ of $r^2 + 2r - C = 0$.
If $x^3z^2 \tends r$ then $|z| = |r|^{1/2}|x|^{-3/2} + o(|x|^{-3/2})$.
Since $5/4 < 3/2$, any
such z automatically satisfies the $z = o(|x|^{-5/4})$ requirement.

To sum up, the following {\it necessary} requirements on the
asymptotic behavior along a $y$-finite asymptotic curve
for $P$ with asymptotic limit $C$ have been derived.
$$ y = x^{-1} + s|x|^{-3/2} + o(|x|^{-3/2}) $$
where $|s|=|r|^{1/2}$ and $r$ is a root of $r^2 + 2r - C = 0$. 
If $r \ne 0$, then only one of the two possible choices of
$s$ occurs for a given asymptotic curve.
However, either choice will lead to a curve with the
right properties, since in either case $x^3z^2 \tends r$.

To verify that these conditions suffice, denote again by
$z$ the (new) error term.
$$ y = x^{-1} + s|x|^{-3/2} + z $$
Compute $P(x,x^{-1} + sx^{-3/2})$. The result is
$$ s^4 + 2s^2 + (3s^3+3s)x^{-1/2} + (3s^2+1)x^{-1} + sx^{-3/2} $$
This is a correct formula for what must happen if $x \tends +\infty$.
If $x \tends -\infty$ instead, $x^{-3/2}$ must be replaced by
$|x|^{-3/2} = (-x)^{-3/2}$. 
Compute $P(x,x^{-1} + s(-x)^{-3/2})$. The result is
$$ s^4 - 2s^2 + (-3s^2+3s)(-x)^{-1/2} + (3s^2-1)(-x)^{-1} + s(-x)^{-3/2} $$
To obtain the corresponding asymptotic identities in
Peretz's standard form, substitute $1/x^2$ for $x$
and $y$ for $s$ in the first,
and $-1/x^2$ for $x$ and $y$ for $s$ in the second,
obtaining
\begin{align}
\notag
P(1/x^2,yx^3+x^2) &= y^4+2y^2+(3y^3+3y)x+(3y^2+1)x^2+yx^3 \\
P(-1/x^2,yx^3-x^2) &= y^4 - 2y^2 + (3y^3-3y)x + (3y^2-1)x^2 + yx^3
\notag
\end{align}

\begin{remark}
In \cite{AsympVals} Peretz claims that to find all the 
asymptotic values of a polynomial $P$ corresponding to $y$-finite asymptotic
curves, it suffices to consider asymptotic identities of the
form $P(1/x^k,yx^N + a_{N-1}x^{N-1} + \ldots + a_0) = a(x,y)
\in \R[x,y]$. This appears to be an oversight. As this case
shows, one must consider asymptotic identities involving
$\pm 1/x^k$ when $k$ is even, otherwise asymptotic values
obtained as $x \tends -\infty$ will be missed. 
It turns out that both $P$ and $Q$ satisfy asymptotic
identities for each of the two asymptotic curves above.
The $(u,v)$
coordinates of points in $\text{AV}[F]$ are obtained by
substituting $x=0$ in the right hand sides of the asymptotic
identities, and allowing $y$ to vary. The right hand side
of the first identity for $P$ reduces to $y^4 + y^2$ for $x=0$,
so one can obtain only points with $u \ge 0$. In fact, one
obtains the points in Figure 1 on the $q+$ portion of the curve,
starting at $(0,0)$ and going to the right.
The remaining points in $\text{AV}[F]$
all derive from the identities for the second
asymptotic curve $(-1/x^2,yx^3-x^2)$.
\end{remark}

Next consider the case in which $z$, the error term, is not
identically zero. 
As an illustration,
compute $P(x,x^{-1} + sx^{-3/2} + z)$.
The result is the same result obtained when $z=0$ plus the
following additional terms
\begin{multline}
 z^4x^6 + 4sz^3x^{9/2} + 3z^3x^4 + (6s^2z^2+2z^2)x^3 \\
+ 9sz^2x^{5/2} + 3z^2x^2 + (4s^3z+4sz)x^{3/2} \\
+ (3z+9s^2z)x + 6szx^{1/2} + z
\notag
\end{multline}
Each of these terms tends to zero as a consequence of $z=o(|x|^{-3/2})$.
So these are indeed all asymptotic curves, and the limiting
asymptotic value obtained is independent of the form of $z$ as long
as $z=o(|x|^{-3/2})$.
No new asymptotic limits are found. However, formally different
asymptotic identities can be derived. For instance, from
$y = x^{-1} + ax^{-3/2} + bx^{-2}$ one obtains 
the following asymptotic identity when $1/x^2$ is
substituted for $x$ and $y$ is substituted for $b$
\begin{multline}
 P(1/x^2,yx^4 + ax^3 + x^2) = \\
x^4y^4 + (4x^3a+3x^4)y^3 + (9x^3a+6a^2x^2+3x^4+2x^2)y^2 \\
+ (6x^3a+x^4+4ax+4a^3x+3x^2+9a^2x^2)y \\
+ a^4 + 3ax + 3a^2x^2 + x^2 + x^3a + 3a^3x + 2a^2
\notag
\end{multline}
Setting $x=0$ in the right hand side to see what asymptotic
limits are obtained yields $a^4 + 2a^2$, the same set of
limit values as for the previous asymptotic identity.
Note that the free parameter $a$ yields the asymptotic
values here, whereas all the $y$ terms disappear if $x$
is set equal to zero.

To look for $x$-finite asymptotic values, similar steps are
taken, but there are fewer such steps since the powers of $y$
extend only up to $y^4$. The first assertion, from the
bottom up, that has a nonzero constant term on the left is
(suitably rearranged)
$$ (x^2y)^3 + (x^2y)^2(-4x+3)+(x^2y)(6x^2-7x+3)+1 \tends 0 $$
and as $x \tends 0$ this implies that $x^2y$ tends to a root $r$
of the equation 
$$r^3 +3r^2 +3r + 1 = (r+1)^3 = 0$$
Thus the first approximation is
$$ x = -y^{-1/2} + o(y^{-1/2})$$
with $y \tends +\infty$ the only possibility.
Substitute $x = -y^{-1/2} + z$ into the Peretz assertions to
obtain
\begin{multline}
z^6y^4+18z^4y^3-4y^3z^5+6z^4y^2-47z^3y^2+36y^2z^2+8y \\
-4z^3y -44zy + 41z^2y + 11 -12z +20z^4y^{(5/2)}+12z^2y^{(1/2)} \\
-24zy^{(3/2)}-34zy^{(1/2)}+61z^2y^{(3/2)}-24z^3y^{(3/2)}-32z^3y^{(5/2)} \\
-6z^5y^{(7/2)}+14y^{(1/2)}+4y^{(-1/2)} \\
\tends 0 \notag \label{b1}  
\end{multline}
\begin{multline}
z^6y^3 + 18z^4y^2 -4z^5y^2 +6z^4y - 47z^3y +36z^2y + 8 \\
+ 36z^2 -41z +6y^{(-1)} + 20z^4y^{(3/2)} -24z^3y^{(1/2)}
-24zy^{(1/2)} \\ + 61z^2y^{(1/2)} -24zy^{(-1/2)} -6z^5y^{(5/2)} -32z^3y^{(3/2)}
+ 11y^{(-1)} \\
\tends 0 \notag \label{b2}  
\end{multline}
\begin{multline}
-18zy^{(-1/2)} + 40z^2y^{(-1/2)} -20zy^{(-1)} + z^6y^2 - 6z^5y^{(3/2)} +18z^4y \\
-32z^3y^{(1/2)} -4yz^5 + 20z^4y^{(1/2)} + 4y^{(-3/2)} - 40z^3
+ 33 z^2 + 4y^{(-1)} \\
\tends 0 \notag \label{b3}  
\end{multline}
\begin{multline}
z^6y - 6z^5y^{(1/2)} + 15z^4 -20z^3y^{(-1/2)} + 15z^2y^{(-1)} -6zy^{(-3/2)} +y^{(-2)} \\
\tends 0 \notag \label{b4}  
\end{multline}
 and the trivial assertion $0 \tends 0$. Every term containing a monomial
of the form $z^my^n$ tends to zero 
if $m \ge 2n$. 
The last two assertions collapse to the trivial $0 \tends 0$.
However, the next one from the bottom up reduces to
$8 \tends 0$. Since that cannot happen, it follows that
there are no $x$-finite asymptotic limits.

\section{The Asymptotic Variety of $F$}
In the previous section it was shown that the asymptotic
curves
$(1/x^2,yx^3+x^2)$ and $(-1/x^2,yx^3-x^2)$, both defined
for $x \ne 0$, are a basis for the asymptotic values of $P$,
in the sense that every asymptotic value of $P$ arises as
a limit along one of these curves. These curves are also
asymptotic curves that yield finite asymptotic values for $Q$.
Specifically, one has the asymptotic identities
\begin{multline}
\notag
Q(1/x^2,yx^3+x^2) =
-(1/4)(1736y^6 + 1044 y^4 + 1155 y^8 + 300 y^10) \\
-(x/4)(5700y^7 + 6692 y^5 + 2792 y^3 + 1800 y^9) \\
-(x^2/4)(9636y^4 + 11250 y^6 + 4500 y^8 + 2432 y^2) \\
-(x^3/4)(6000y^7 + 680y + 11100y^5 + 6140 y^3) \\
-(x^4/4)(4500y^6 + 1460y^2 + 5475y^4) \\
-(x^5/4)(1800y^5 + 1080y^3) -75x^6y^4
\end{multline}
\begin{multline}
\notag
Q(-1/x^2,yx^3-x^2) =
+(1/4)(1736y^6 - 1044 y^4 - 1155 y^8 + 300 y^10) \\
+(x/4)(-5700y^7 + 6692 y^5 - 2792 y^3 + 1800 y^9) \\
+(x^2/4)(9636y^4 - 11250 y^6 + 4500 y^8 - 2432 y^2) \\
+(x^3/4)(6000y^7 - 680y - 11100y^5 + 6140 y^3) \\
+(x^4/4)(4500y^6 + 1460y^2 - 5475y^4) \\
+(x^5/4)(1800y^5 - 1080y^3) + 75x^6y^4
\end{multline}
Substituting $x=0$ to obtain the asymptotic values,
both here and in the asymptotic identities for $P$,
yields the following two parameterized curves that
together make up the whole asymptotic variety
$$ u = y^4 + 2y^2,\; v = -(1/4)(1736y^6 + 1044 y^4 + 1155 y^8 + 300 y^{10}) $$
$$ u = y^4 - 2y^2,\; v =  (1/4)(1736y^6 - 1044 y^4 - 1155 y^8 + 300 y^{10}) $$
These two parameterizations can be combined into one by
putting $y^2=h$ (for $h \ge 0$) in the first,
and $y^2 = -h$ (for $h \le 0$) in the second,
which yields exactly the parameterization
considered before 
$$u = h^2 + 2h,\; v = -(1/4)(1736h^3 + 1044 h^2 + 1155 h^4 + 300 h^5) $$

\begin{remark}
The functions $t$,$h$, and $f$
introduced in the definition of $P$ and
$Q$ all satisfy asymptotic identities with respect to each of the above
two asymptotic curves. As Peretz remarked in \cite{AsympVals,OnCxKeller},
examples such as Pinchuk's arise from finding pairs of polynomials with
a nowhere vanishing Jacobian determinant in a real subalgebra of $\R[x,y]$
consisting of polynomials all of which have one or more shared asymptotic
curves for which they satisfy asymptotic identities.
\end{remark}

\section{The Asymptotic Flower of $F$}
In \cite{Flower} the authors consider (primarily polynomial) maps
of the real plane to itself that are proper. The {\it flower} of
a map is the inverse image of the set of critical values (and a value
is critical precisely if it is the image of a point at which the
Jacobian determinant vanishes). Away from the flower the map is locally
a covering map (proper and a local homeomorphism). In fact, it is a
a covering map 
(over its image)
on any connected component of the complement of the flower.
In the case of Pinchuk's map the flower as defined above is empty, but
the covering property fails to hold because the map is not proper.
This suggests calling the above flower the {\it critical flower} and
introducing as well the {\it asymptotic flower}, defined as the inverse
image of the set of asymptotic values of the map. On the complement of
the asymptotic flower, the restricted map to the complement of the set
of asymptotic values is proper. 
This is because, by definition, as the asymptotic flower is approached,
the image of a point will tend to infinity in the image plane or to an
asymptotic value.
But since the codomain of the restricted map is the complement of
the set of asymptotic values, this means that the image is tending
to infinity relative to that codomain.
(Note. Asymptotic values can be defined
in terms of limits of sequences as well on manifolds, using local pathwise
connectedness to produce the appropriate curves. More general definitions
are possible as well.) On each component of the complement of the
{\it total flower} (the union of the critical and asymptotic flowers) the
restricted map to the complement of the critical and asymptotic values will
be a covering map over its image.

\begin{figure}[hbt]
\vspace{-0.5in}
\centerline{
\psfig{figure=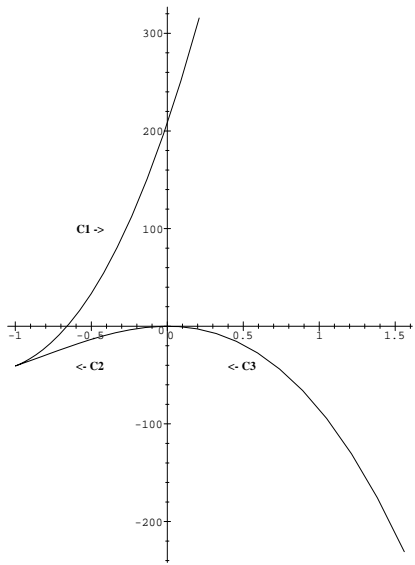,height=4in,width=3in}
}
\caption{Three curves in the asymptotic variety.}
\label{AVmorefig}
\end{figure}

Since the asymptotic variety $\text{AV[F]}$ has been identified for
the $F$ at hand, it remains only to compute its inverse image to
obtain $\text{AF}[F]$, the asymptotic flower of $F$. By consulting Table 2,
it becomes clear that
\begin{itemize}
\item exactly two points, $(0,0)$ and $(-1,-163/4)$,
have no inverse images
\item every other point of $\text{AV}[F]$ has exactly one inverse image
\item every point not in $\text{AV}[F]$ has exactly two inverse images
\end{itemize}
This follows from a case by case check, cases corresponding to rows of
the table.
If one removes the two points that have no inverse images from
$\text{AV}[F]$, it breaks up into three connected curves. Call
them $C1,C2,C3$, as follows.

 $C1$ is the $q-$ curve, starting
at $(-1,-163/4)$ and continuing up and to the right. $C2$ is the
portion of the $q+$ curve starting at $(0,0)$ and continuing down
and to the left, ending at $(-1,-163/4)$. Finally, $C3$ is the
portion of the $q+$ curve ending at $(0,0)$ and arriving from down
and to the right. Starting and ending points mentioned are not
actually points of these curves, since they represent precisely
points that were removed.
The descriptions also imply orientations for the the three curves.
Figure 2 shows the curves and their orientations.

Each point of each of the three curves has exactly one inverse image.
Furthermore, as a starting or ending point, finite or infinite, is
approached, the inverse image point tends to infinity. Thus the inverse
image of each of $C1,C2,C3$ is a smooth curve in the plane (no singularities
and no self-intersections) that tends to infinity at either end.
Call these curves $D1,D2,D3$. 
By definition, $\text{AF}[F] = D1 \cup D2 \cup D3$, where the
curves are considered as point sets.

\begin{figure}[hbt]
\vspace{-0.5in}
\centerline{
\psfig{figure=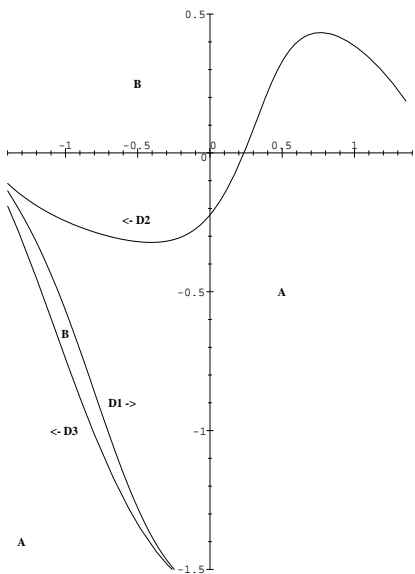,height=4in,width=3in}
}
\caption{The asymptotic flower of Pinchuk's map.}
\label{AFfig}
\end{figure}

Each of these curves $D1,D2,D3$ divides the plane into
two simply connected parts (the Jordan curve Theorem), which may be described
as the regions left and right of the curve, using the induced orientations
to define left and right. Removing the curves $D1,D2,D3$ thus leaves exactly
four simply connected open components. The restriction $F: \R^2 \setminus
\text{AF}[F] \to \R^2 \setminus \text{AV}[F]$, maps each component into
either the region $L$ to the left of $\text{AV}[F]$ 
or into the region $R$ to its right. 
Each region mapped
into $L$ is, in fact, mapped homeomorphically onto $L$, because we are
dealing with a covering of a simply connected region. 
Similarly for $R$.
Label a connected component of $\R^2 \setminus \text{AF}[F] =
\R^2 \setminus (D1 \cup D2 \cup D3)$ with $A$ if it maps to
$L$, and with $B$ if it maps to $R$.
Two of the regions
are labeled $A$, and two are labeled $B$. 
The global data of the map (the relations between the domains
and curves) are best explained by a figure. Figure 3 depicts
the component curves of $\text{AF}[F]$ and their orientations,
and also labels the regions defined by the curves. Figure 3 uses
nonlinear scaling to produce a more comprehensible picture; the
values plotted are actually the arctangents of the coordinates
$x$ and $y$. The figure was generated by solving for
the inverse images of a large number of points on
$\text{AV}[F]$.
The labeling of the regions can be checked by computing
the images of a few points not in the flower (and can
also be deduced to a large extent from the fact that $F$
is orientation preserving, since its Jacobian determinant
is everywhere positive).

\bibliographystyle{amsplain}
\bibliography{jc,samuelson,algeo,diffeq}

\ifx\undefined\bysame
\newcommand{\bysame}{\leavevmode\hbox to3em{\hrulefill}\,}
\fi
\begin{thebibliography}{1}

\bibitem{PPR2}
L.~Andrew Campbell, {\em Partial properness and real planar maps}, Applied
  Math. Letters {\bf 9} (1996), no.~5, 99--105.

\bibitem{Flower}
Iaci Malta, Nicolau~C. Saldanha, and Carlos Tomei, {\em The numerical inversion
  of functions from the plane to the plane}, Mathematics of Computation {\bf
  65} (1996), no.~216, 1531--1552.

\bibitem{OnCxKeller}
Ronan Peretz, {\em On counterexamples to {K}eller's problem}, Illinois J. Math.
  {\bf 40} (1996), no.~2, 293--303.

\bibitem{AsympVals}
\bysame, {\em The variety of asymptotic values of a real polynomial etale map},
  Journal of Pure and Applied Algebra {\bf 106} (1996), 102--112.

\bibitem{StrongRealJC}
Sergey Pinchuk, {\em A counterexample to the strong real {J}acobian
  conjecture}, Math. Zeitschrift {\bf 217} (1994), 1--4.

\bibitem{Pinchuk.Cx.Email}
Arno van~den Essen, {\em Personal communication}, Email, June 1994, Provided
  details of a degree (10,25) case of Pinchuk's counterexample to the Real
  Jacobian Conjecture.

\end{thebibliography}
\end{document}